\documentclass{elsarticle}
\usepackage[ansinew]{inputenc}
\usepackage{amsfonts}
\usepackage{latexsym}
\usepackage{amsmath}
\usepackage{amssymb}

\newcommand{\supp}{\textrm{supp }}

\newcommand{\C}{\mathbb{C}}
\newcommand{\N}{\mathbb{N}}

\renewcommand{\P}{\mathbb{P}}

\newtheorem{defin}{Definition}[section]

\newtheorem{theorem}[defin]{Theorem}
\newtheorem{exa}[defin]{Example}

\newtheorem{lemma}[defin]{Lemma}

\newtheorem{corollary}[defin]{Corollary}
\newtheorem{qu}{Question}

\newenvironment{proof}
{\noindent{\it Proof.\,}}{\hfill $\Box$\par\vspace{2.5mm}}

\numberwithin{equation}{section}

\begin{document}
\begin{frontmatter}

\title{Difference Picard theorem for meromorphic functions of several variables\tnoteref{t1}}

\tnotetext[t1]{The research reported in this paper was supported
in part by the Academy of Finland grant \#118314 and \#210245, the Isaac Newton Institute for Mathematical Sciences, and
the NordForsk foundation.}

\author{Risto Korhonen}
\ead{risto.korhonen@helsinki.fi}
\address{Department of Mathematics and Statistics, P.O.~Box~68, FI-00014
University of Helsinki, Finland}


\begin{abstract}
It is shown that if $n\in\N$, $c\in\C^n$, and three distinct
values of a meromorphic function $f:\C^n\to\P^1$ of hyper-order $\varsigma(f)$
strictly less than $2/3$ have forward invariant pre-images with
respect to a translation $\tau:\C^n\to\C^n$, $\tau(z)=z+c$, then
$f$ is a periodic function with period~$c$. This result can be
seen as a generalization of M. Green's Picard-type theorem in the special case where $\varsigma(f)<2/3$,
since the empty pre-images of the usual Picard exceptional
values are by definition always forward invariant. In
addition, difference analogues of the lemma on the logarithmic
derivative and of the second main theorem of Nevanlinna theory for
meromorphic functions $\C^n\to\P^1$ are given, and their
applications to partial difference equations are discussed.
\end{abstract}


\end{frontmatter}

\section{Introduction}

The purpose of this paper is to find difference analogues of the
lemma on the logarithmic derivative and of the second main theorem
of Nevanlinna theory for meromorphic functions, where the
operation of partial differentiation in the ramification term has
been replaced by the genuine shift operator $\Delta_c
f:=f(z_1+c_1,\ldots,z_n+c_n)-f(z_1,\ldots,z_n)$,
$c=(c_1,\ldots,c_n)\in\C^n$, operating on a meromorphic function
$f:\C^n\to\P^1$ of hyper-order strictly less than $2/3$.
Hyper-order is defined by
    \begin{equation}\label{hyper}
    \varsigma(f)=\limsup_{r\to\infty}\frac{\log\log T(r,f)}{\log r},
    \end{equation}
where $T(r,f)$ is the Nevanlinna characteristic function of $f$
(see Section~\ref{results} below for a short review of Nevanlinna
theory of several variables). These results will have two main
applications. First, we will obtain a difference analogue of
Picard's theorem in several variables, which says that if
$n\in\N$, $c\in\C^n$, and three distinct values of a meromorphic
function $f:\C^n\to\P^1$ such that $\varsigma(f)<2/3$ have forward
invariant pre-images with respect to a translation
$\tau:\C^n\to\C^n$, $\tau(z)=z+c$, then $f$ is a periodic function
with period~$c$. In the special case of $\varsigma(f)<2/3$ this
result can be seen as a generalization of M. Green's Picard-type
theorem, since the (empty) pre-images of the usual Picard
exceptional values are special cases of forward invariant
pre-images. The second application can be described as a Malmquist
type theorem for partial difference equations. We will show that
the existence of one meromorphic solution $w:\C^n\to\P^1$ such
that $\varsigma(w)<2/3$ is enough to reduce a large class of
partial difference equations into a difference Riccati equation.

The remainder of the paper is organized as follows. The difference
analogue of Picard's theorem (Theorem \ref{picard} below) is
stated in Section~\ref{picardsec}. Section \ref{results} contains
difference analogues of the lemma on the logarithmic derivative
and of the second main theorem (Theorems \ref{logdiff} and
\ref{2mthm} below). Applications of these results to partial
difference equations are discussed in Section~\ref{partialsec}.
The difference analogue of the lemma on the logarithmic derivative
in several variables is proved in Sections \ref{estimsec} and
\ref{proofsec}, while Section \ref{proof2nd} contains the proof of
the difference version of the second main theorem. Finally, the
difference analogue of Picard's theorem in several variables is
proved in Section \ref{picardproof}.

\section{A difference analogue of Picard's theorem}\label{picardsec}

Picard's theorem states that any non-constant entire function
$f(z)$ assumes all values in the complex plane with at most one
possible exception \cite{picard:79}. Fatou \cite{fatou:22b,fatou:22} has constructed an example of a
biholomorphic mapping $f:\C^2\to\C^2$ such that the set of Picard
exceptional values $\C^2\setminus f(\C^2)$ contains a non-empty
open set. At first sight this example appears to imply severe difficulties in
generalizing Picard's theorem to meromorphic functions of several
variables. However, it turns out that there is a natural generalization which can be found
by rephrasing Picard's theorem in terms of projective spaces.
Green \cite{green:72} showed that any holomorphic mapping from
$\C^n$ into the projective space $\P^m$ that misses $2m+1$
hyperplanes in general position is a constant, thus improving an
earlier Picard-type theorem by Wu \cite{wu:69}. Moreover, extensions of
Nevanlinna's second main theorem to several variables can be
regarded as deep generalizations of Picard's theorem, see, for
instance,
\cite{carlsong:72,griffithsk:73,wongs:94,ye:95,cherryy:97}.

We will show that forward invariance with respect to a translation
of the pre-image of a target value is, in the sense of Picard
exceptionality, as restrictive for non-periodic meromorphic
functions $\C^n\to\P^1$ such that $\varsigma(f)<2/3$, as omitting
the target value completely. We say that the pre-image of
$a\in\P^1$ is under $f$ is forward invariant with respect to the
translation $\tau$ if $\tau(f^{-1}(\{a\}))\subset f^{-1}(\{a\})$
where $\tau(f^{-1}(\{a\}))$ and $f^{-1}(\{a\})$ are considered to
be multisets in which each point is repeated according to its
multiplicity. By this definition the (empty and thus forward
invariant) pre-images of the usual Picard exceptional values
become special cases of forward invariant pre-images. The
following theorem is a difference analogue of Picard's theorem for
meromorphic functions in several variables.

\begin{theorem}\label{picard}
Let $f:\C^n\to\P^1$ be a meromorphic function such that $\varsigma(f)<2/3$, and
let $\tau(z)=z+c$, where $\tau:\C^n\to\C^n$ and $c\in\C^n$.  If three distinct values of $f$ have
forward invariant pre-images with respect to $\tau$, then $f$ is a periodic function with
period~$c$.
\end{theorem}

Theorem \ref{picard} is proved in Section~\ref{picardproof} below.
A simple example from \cite{halburdkt:09} shows that the condition on growth of $f$ cannot
be removed, at least not completely. By taking
$g(z)=\exp(\exp(z))$, the pre-image of each of the
$m^{\textrm{th}}$ roots of unity is forward invariant with respect
to the translation $\tau(z)= z+\log(m+1)$. Since clearly
$g(z)\not\equiv g(z+\log(m+1))$, it follows that a slightly weaker
growth condition in Theorem \ref{picard} would allow a
non-periodic meromorphic function with \textit{arbitrarily many}
values having forward invariant pre-images.

\section{Second main theorem}\label{results}

%

One of the key components in Nevanlinna's original proof of the
second main theorem is a technical result referred to as
the lemma on the logarithmic derivative. This lemma has also been
used as an important tool in the study of value distribution of
meromorphic solutions of differential equations in the complex
plane \cite{hayman:64,laine:93,gromakls:02}. The original proof of
the second main theorem in several variables was based on a
differential geometric method due to Ahlfors and F. Nevanlinna,
see, e.g., \cite{wongs:94}, instead of Nevanlinna's method based on the
lemma on the logarithmic derivative. The first generalization of
the lemma on the logarithmic derivative to several complex
variables was given by Vitter \cite{vitter:77}, who used the
method of non-negative curvature developed by Carlson, Cowen, Griffiths and King \cite{carlsong:72,griffithsk:73,coweng:76}.
Biancofiore and Stoll used an alternative method based
on a technique they call ``fiber integration'' to prove their version of the lemma on the
logarithmic derivative in several complex variables
\cite{biancofiores:81}. Further improvements and generalizations
of the lemma on the logarithmic derivative has been given, for
instance, by Cherry \cite{cherry:93} and Ye \cite{ye:95,ye:96}.

The purpose of this section is to present difference analogues of
the lemma on the logarithmic derivative and of the second main
theorem in several complex variables. Before stating these two key
results of this paper, we will briefly recall some of the standard
notation of Nevanlinna theory in $\C^n$
\cite{stoll:74,kujalav:74,stoll:85} (see also, for instance,
\cite{vitter:77,biancofiores:81,ye:96}).

Let $z=(z_1,\ldots,z_n)\in\C^n$, and let $r>0$. Introducing the
differential operators $d:=\partial+\overline{\partial}$ and
$d^c:=(\partial-\overline{\partial})/4\pi i$, we define
$\omega_n(z):=dd^c\log|z|^2$ and
$\sigma_n(z):=d^c\log|z|^2\wedge\omega_n^{n-1}(z)$ where
$z\in\C^n\setminus{\{0\}}$ and $|z|^2:=|z_1|^2+\cdots+|z_n|^2$.
Then $\sigma_n(z)$ defines a positive measure with total measure
one on the boundary $\partial B_n(r):=\{z\in\C^n:|z|=r\}$ of the
ball $B_n(r):=\{z\in\C^n:|z|<r\}$. In addition, by defining
$\upsilon_n(z):=dd^c|z|^2$ and $\rho_n(z):=\upsilon_n^n(z)$ for
all $z\in\C^n$, it follows that $\rho_n(z)$ is the Lebesgue
measure on $\C^n$ normalized such that the ball $B_n(r)$ has
measure $r^{2n}$.

Let $f$ be a meromorphic function in $\C^n$ in the sense that $f$
can be written as a quotient of two relatively prime holomorphic
functions. We will write $f=(f_0,f_1)$ where $f_0\not\equiv0$, and
regard $f$ as a meromorphic map $f:\C^n\to\P^1$ such that
$f^{-1}(\infty)\not=\C^n$. The standard definition of Nevanlinna
characteristic function of $f$ is given by
\begin{equation*}
T_f(r,s):=\int_s^r\frac{A_f(t)}{t}\,dt
\end{equation*}
where $0<s<r$ and
\begin{equation*}
A_f(t)=\frac{1}{t^{2n-2}}\int\limits_{B_n(t)}f^*\omega
\wedge\upsilon_n^{n-1} =\int\limits_{B_n(t)}f^*\omega \wedge
\omega_n^{n-1}+A_f(0)
\end{equation*}
is a measure of the spherical area covered by the image of $B_n(t)$ under $f$. Here the pullback
$f^*\omega$ satisfies
\begin{equation*}
f^*\omega = dd^c\log(|f_0|^2+|f_1|^2)
\end{equation*}
for all $z$ outside of the set of indeterminacy
$I_f:=\{z\in\C^n:f_0(z)=f_1(z)=0\}$ of $f$.

A \textit{divisor} on $\C^n$ is an integer valued function which
is locally the difference between the zero-multiplicity functions
of two holomorphic functions, in our case $f_0$ and $f_1$. Let
$a\in\P^1$ such that $f^{-1}(a)\not=\C^n$. Then the $a$-divisor
$\nu_f^a$ of $f=(f_0,f_1)$ is the divisor associated to the
holomorphic functions $f_1-af_0$ and $f_0$. By denoting
$S(r):=\overline{B}_n(r)\cap \supp \nu_f^a$, where
$\overline{B}_n(r)=\{z\in\C^n:|z|\leq r\}$ and $\supp \nu_f^a$
denotes the closure of the set $\{z\in\C^n:\nu_f^a(z)\not=0\}$, we
may define the \textit{counting function} of $\nu_f^a$ as
\begin{equation*}
n_f(r,a):=r^{2-2n}\int\limits_{S(r)} \nu_f^a \upsilon_n^{n-1}
\end{equation*}
for all $n\geq1$ and for all $r>0$.

There are slightly different ways to continue the formulation of
Nevanlinna theory from here. Stoll \cite{stoll:85} defines the
(integrated) \textit{counting function} of $\nu_f^a$ as
\begin{equation*}
N_f(r,s,a):=\int_s^r \frac{n_f(t,a)}{t}\,dt
\end{equation*}
for all $0<s<r$, and the \textit{compensation function} as
\begin{equation*}
m_f(r,a):=\int\limits_{\partial B_n(r)}\hspace{-2mm}
\log\frac{1}{\|f,a\|}\,\sigma_n(z),
\end{equation*}
where $||f,a||$ denotes the chordal distance from $f$ to
$a\in\P^1$. Then the first main theorem of Nevanlinna theory
becomes
\begin{equation*}
T_f(r,s) = N_f(r,s,a) + m_f(r,a) - m_f(s,a)
\end{equation*}
where $0<s<r$.

We choose a slightly different approach (see e.g. \cite{ye:96}) by
denoting $N(r,f):=N_f(r,0,\infty)$ and $N(r,1/(f-a)):=N_f(r,0,a)$,
where $a\not=\infty$ and we have assumed that $f(0)\not=a,\infty$.
Then by the Jensen formula,
\begin{equation}\label{Jensen}
N\left(r,\frac{1}{f}\right)-N(r,f)=\int\limits_{\partial B_{n}(r)} \log |f(z)|\sigma_n(z)-\log|f(0)|
\end{equation}
for all $r>0$, provided that $f(0)\not=0,\infty$. By defining the \textit{proximity function} of $f$ as
\begin{equation*}
m(r,f):=\int\limits_{\partial B_n(r)}\hspace{-2mm}
\log^{+}|f(z)|\sigma_n(z),
\end{equation*}
and if $a\not=\infty$,
\begin{equation*}
m\left(r,\frac{1}{f-a}\right):=\int\limits_{\partial
B_n(r)}\hspace{-2mm} \log^{+}\frac{1}{ |f(z)-a|}\sigma_n(z),
\end{equation*}
the Jensen formula \eqref{Jensen} becomes
\begin{equation}\label{1mt}
T(r,f)=m\left(r,\frac{1}{f-a}\right)+N\left(r,\frac{1}{f-a}\right)-\log\frac{1}{|f(0)-a|}
\end{equation}
where $T(r,f)=m(r,f)+N(r,f)$ and $f$ is a meromorphic function on
$\C^n$ satisfying $f(0)\not=a,\infty$. The order of growth of $f$
is defined by
\begin{equation*}
\rho(f):=\limsup_{r\to\infty}\frac{\log T(r,f)}{\log r}.
\end{equation*}
The following theorem is a difference analogue of the lemma on the
logarithmic derivative in several complex variables. It
generalizes the one dimensional result \cite[Theorem
2.1]{halburdk:06AASFM} by Halburd and the author. Recall the
definition of hyper-order from \eqref{hyper}.

\begin{theorem}\label{logdiff}
Let $f$ be a non-constant meromorphic function in
$\C^n$ such that $f(0)\not=0,\infty$, let $c\in\C^n$, and let
$\varepsilon>0$. If $\varsigma(f)=\varsigma<2/3$, then
\begin{equation}\label{mainrel}
\int\limits_{\partial B_n(r)}
\log^{+}\left|\frac{f(z+c)}{f(z)}\right|\sigma_n(z) =
o\left(\frac{T(r,f)}{r^{1-\frac{3}{2}\varsigma-\varepsilon}}\right)
\end{equation}
for all $r>0$ outside of a possible exceptional set
$E\subset[1,\infty)$ of finite logarithmic measure $\int_E
1/dt<\infty$.
\end{theorem}

The proof of Theorem \ref{logdiff} can be found in Sections
\ref{estimsec} and \ref{proofsec} below. Recall that we have
adopted the notation $\Delta_c f:=f(z+c)-f(z)$ for $c\in\C^n$ and
$f:\C^n\to\P^1$. The following corollary is an immediate
consequence of Theorem~\ref{logdiff}.

\begin{corollary}\label{shiftcor}
Let $a$ and $c$ be constants in $\C^n$, let $f$ be a non-constant
meromorphic function in $\C^n$ such that
$f(0)\not=a,\infty$, and let $\varepsilon>0$. If $\varsigma(f)=\varsigma<2/3$, then
\begin{equation*}
m\left(r,\frac{\Delta_c f}{f-a}\right)=
o\left(\frac{T(r,f)}{r^{1-\frac{3}{2}\varsigma-\varepsilon}}\right)
\end{equation*}
for all $r>0$ outside of a possible exceptional set
$E\subset[1,\infty)$ of finite logarithmic measure.
\end{corollary}

Corollary \ref{shiftcor} can be applied to prove a difference
analogue of the second main theorem of Nevanlinna theory for
meromorphic functions $f:\C^n\to\P^1$, which extends \cite[Theorem
2.4]{halburdk:06AASFM} to meromorphic functions of several
variables.


\begin{theorem}\label{2mthm}
Let $c\in\C^n$, let $\varepsilon>0$, and let $f$ be a meromorphic
function in $\C^n$ such that $\Delta_c
f\not\equiv0$. Let $q\geq 2$, and let $a_1,\ldots,a_q\in \P^1$ be
distinct finite constants such that $f(0)\not=a_j,\infty$ for all
$j=1,\ldots,q$. If $\varsigma(f)=\varsigma<2/3$, then
\begin{equation*}
m(r,f)+\sum_{j=1}^q m\left(r,\frac{1}{f-a_j}\right)\leq
2T(r,f)-N_{\Delta}(r,f)+o\left(\frac{T(r,f)}{r^{1-\frac{3}{2}\varsigma-\varepsilon}}\right),
\end{equation*}
where
\begin{equation*}
N_{\Delta}(r,f)=2N(r,f)-N(r,\Delta_c f)+N\left(r,\frac{1}{\Delta_c f}\right),
\end{equation*}
and $r$ lies outside of a possible exceptional set
$E\subset[1,\infty)$ of finite logarithmic measure.
\end{theorem}

The proof of Theorem \ref{2mthm} can be found from Section
\ref{proof2nd}.

\section{Applications to partial difference equations}\label{partialsec}

Ablowitz, Halburd and Herbst \cite{ablowitzhh:00} have suggested
that the existence of sufficiently many finite-order meromorphic
solutions could be used as a detector of Painlev\'e type
difference equations. Halburd and the author used one-dimensional
difference analogues \cite{halburdk:06JMAA,halburdk:06AASFM} of
some of the main results of Nevanlinna theory to prove that the
existence of at least one finite-order meromorphic solution, which
is not simultaneously a solution of a first-order difference
Riccati equation, is enough reduce a large class of difference
equations into a list of equations consisting exactly of known
discrete equations of Painlev\'e type
\cite{halburdk:07PLMS,halburdk:07JPA}.




The purpose of this section is to extend some of the methods used
in \cite{ablowitzhh:00} to partial differences, and apply these
generalized results to single out the difference Riccati equation
out of a large class of first-order partial difference equations.
We start by stating the main result of this section.

Let $\mathcal{S}(f)=\{g:\C^n\to\P^1\textrm{ meromorphic
}:T(r,g)=o(T(r,f))\}$ where $r\to\infty$ outside of a possible
exceptional set of finite logarithmic measure. A meromorphic
solution $w:\C^n\to\P^1$ of a difference equation is called
\textit{admissible} if all coefficients of the equation are in
$\mathcal{S}(f)$ (see \cite[p.~192]{laine:93}).

\begin{theorem}\label{riccati}
Let $c\in\C^n$. If the difference equation
    \begin{equation}\label{wR}
    w(z+c) = R(z,w(z)),
    \end{equation}
where $R(z,u)$ is rational in $u$ having meromorphic coefficients
in $\C^n$, has an admissible meromorphic solution $w:\C^n\to\P^1$
such that $\varsigma(w)<2/3$, then $\deg_w(R)=1$.
\end{theorem}


The first result needed in the proof of Theorem \ref{riccati} is
due to Valiron \cite{valiron:31} and Mohon'ko \cite{mohonko:71}.

\begin{theorem}[\cite{mohonko:71,valiron:31}]\label{vm}
Let $R(z,u)$ be a rational function of $u$ whose coefficients are
meromorphic functions $h(z)$ in $\C^n$ satisfying
$T(r,h)=O(\phi(r))$, where $\phi$ is a fixed positive increasing
function on $[0,\infty)$. Then for every meromorphic function $f:\C^n\to\P^1$
we have
    \begin{equation*}
    T(r,R(z,f(z)))=\deg_f T(r,f)+O(\phi(r)).
    \end{equation*}
\end{theorem}

According to an identity due to Valiron \cite{valiron:31} and
Mohon'ko \cite{mohonko:71} (see also, e.g.,
\cite[p.~31]{goldbergo:08} and \cite[p.~29]{laine:93})
    \begin{equation}\label{vm2}
    \deg_f(R)T(r,f)=T(r,R(z,f(z)))+O(\phi(r)),
    \end{equation}
whenever $f$ is a non-constant meromorphic function in the complex
plane. As was observed in \cite[Appendix B., p.
453]{goldbergo:08}, by following the proof of \eqref{vm2} in
\cite{mohonko:71} (see also \cite{laine:93}) it can be seen that
the identity \eqref{vm2} holds for any non-decreasing
characteristic function $T(r,f)$ which satisfies the basic
Nevanlinna inequalities, the first main theorem, and the property
$T(r,f^2)=2T(r,f)$. Therefore, in particular, the assertion of Theorem~\ref{vm} follows.

Chiang and Feng \cite{chiangf:08} have shown that if $f$ is a
finite-order meromorphic function in the complex plane and
$\eta\in\C$, then
    \begin{equation}\label{cf}
    T(r,f(z+\eta))=T(r,f)+O(r^{\rho-1+\varepsilon}),\quad
    r\to\infty,
    \end{equation}
where $\rho=\rho(f)$ is the order of $f$ and $\varepsilon>0$. A
similar estimate
    \begin{equation}\label{hk}
    T(r,f(z+\eta))=T(r,f)+o(T(r,f)),
    \end{equation}
where $r\to\infty$ outside of an exceptional set of finite
logarithmic measure, follows by combining \cite[Theorem
2.1]{halburdk:06AASFM} with \cite[Lemma 2.1]{halburdk:07PLMS}. The
following theorem is a generalization of the asymptotic relations
\eqref{cf} and \eqref{hk} to several variables.

\begin{theorem}\label{Trfc}
Let $f:\C^n\to\P^1$ be a meromorphic function, let $c\in\C^n$ and
let $\varepsilon>0$. If $\varsigma(f)=\varsigma<2/3$, then
    \begin{equation}\label{Trfceq}
    T(r,f(z+c))=T(r,f)+o\left(\frac{T(r,f)}{r^{1-\frac{3}{2}\varsigma-\varepsilon}}\right)
    \end{equation}
where $r\to\infty$ outside of an exceptional set of finite
logarithmic measure.
\end{theorem}

\begin{proof}
First we observe that $N(r,f(z+c))\leq N(r+|c|,f)$ by the
definition of the counting function. Therefore, by defining
    $$
    \lambda_2:=\limsup_{r\to \infty}\frac{\log\log N(r,f)}{\log r}
    $$
and applying \cite[Lemma~8.3]{halburdkt:09}, it follows that
    \begin{equation}\label{Nfc}
    N(r,f(z+c))\leq
    N(r+|c|,f)=N(r,f)+o\left(\frac{N(r,f)}{r^{1-\lambda_2-\varepsilon}}\right),
    \end{equation}
where $r$ tends to infinity outside of an exceptional set of
finite logarithmic measure. Second, by Theorem~\ref{logdiff} we
have
    \begin{equation}\label{mfc}
     m\left(r,\frac{f(z+c)}{f(z)}\right) = o\left(\frac{T(r,f)}{r^{1-\frac{3}{2}\varsigma-\varepsilon}}\right)
    \end{equation}
where $r$ lies again outside of an exceptional set of finite
logarithmic measure. The upper bound in the asymptotic relation
\eqref{Trfceq} follows by combining \eqref{Nfc} and \eqref{mfc}
with the inequality
    \begin{equation*}
    T(r,f(z+c))\leq  N(r+|c|,f) + m(r,f) +
    m\left(r,\frac{f(z+c)}{f(z)}\right),
    \end{equation*}
and using the facts $\lambda_2\leq\varsigma$ and $N(r,f)\leq
T(r,f)$. The lower bound follows similarly by combining
    \begin{equation*}
    T(r,f(z))\leq  N(r+|c|,f(z+c)) + m(r,f(z+c)) +
    m\left(r,\frac{f(z)}{f(z+c)}\right),
    \end{equation*}
with \eqref{Nfc} and \eqref{mfc}, applied with the function
$f(z+c)$ and the shift~$-c$.
\end{proof}



We are now ready to prove Theorem~\ref{riccati}.

\medskip

\textit{Proof of Theorem \ref{riccati}. } Suppose that \eqref{wR}
has a meromorphic solution $w:\C^n\to\P^1$ such that
$\varsigma(w)<2/3$. By applying Theorems~\ref{vm} and \ref{Trfc}
respectively to the right and left sides of \eqref{wR}, it follows
that
    \begin{equation*}
    T(r,w)=\deg_w(R)T(r,w)+o(T(r,w))
    \end{equation*}
as $r\to\infty$ outside of an exceptional set of finite logarithmic measure. Therefore, $\deg_w(R)=1$. \hfill $\Box$\par\vspace{2.5mm}

\section{Estimates on integrated difference quotients in $\C$ and $\C^n$}\label{estimsec}

In this section we lay the foundations for the proof of Theorem~\ref{logdiff} by obtaining growth estimates for
integrated difference quotients of a meromorphic function $f$ in $\C$ and in $\C^n$. We start with the one-dimensional case.

\begin{lemma}\label{details1}
Let $f(z)$ be a meromorphic function in $\C$ such that
$f(0)\not=0,\infty$, and let $c\in\C$ and $\delta\in(0,1)$. Then
for all $r>0$ and $s>r+|c|$,
    \begin{equation*}
    \begin{split}
     &\int\limits_{\partial B_1(r)}\hspace{-2mm}\log^{+}\left|\frac{f(z+c)}{f(z)}
    \right|\sigma_1(z) \leq \frac{8\pi|c|^\delta}{\delta(1-\delta)r^\delta}
    \left(n(s,f)+n\left(s,\frac{1}{f}\right)\right) \\
    &\quad +
    \frac{4\pi|c|}{(1-\delta)(s-r-|c|)}\cdot
    \left(\frac{s}{s-r}\right)^{1-\delta}\left(m(s,f)+m\left(s,\frac{1}{f}\right)\right).
    \end{split}
    \end{equation*}
\end{lemma}

Similar estimates to Lemma~\ref{details1} have been obtained
before in \cite[Lemma~2.3]{halburdk:06JMAA}, \cite[Theorem
2.4]{chiangf:08} and \cite[Lemma~8.2]{halburdkt:09} by using
similar methods to here. The improved factor in front of the
function $m(s,f)+m(s,1/f)$ in Lemma~\ref{details1} enables us to
get the inequality \eqref{Im} below in the proof of
Theorem~\ref{logdiff}, instead of a weaker estimate which would
follow by using, for instance, \cite[Lemma~8.2]{halburdkt:09}. The
reason why this is important is the fact that the estimate
\eqref{Im} is ultimately the cause for the slightly unsatisfactory
growth condition $\varsigma(f)<2/3$ in Theorem~\ref{logdiff}. By
applying \cite[Lemma~8.2]{halburdkt:09} instead of
Lemma~\ref{details1} we would be lead to the condition
$\varsigma(f)<2/5$. This also means that if one is interested in
extending Theorem~\ref{logdiff} to meromorpfic functions of
hyper-order less than one, say, then inequality \eqref{Im} is a good place to start looking for potential improvements.

\medskip

\textit{Proof of Lemma \ref{details1}. } The Poisson-Jensen
formula \cite[Theorem 1.1]{hayman:64} implies
    \begin{equation}\label{integratethis}
    \begin{split}
    \log \left|\frac{f(z+c)}{f(z)}\right| &= \int_0^{2\pi}
    \log|f(se^{i\theta})|\textrm{Re}\left(\frac{se^{i\theta}+z+c}{se^{i\theta}-z-c}-
    \frac{se^{i\theta}+z}{se^{i\theta}-z}\right)\,\frac{d\theta}{2\pi}\\
    &\quad  + \sum_{|a_n|<s} \log \left|\frac{s(z+c-a_n)}{s^2-\bar a_n(z+c)}\cdot\frac{s^2-
    \bar    a_n z}{s(z-a_n)}\right| \\
    &\quad  - \sum_{|b_m|<s} \log \left|\frac{s(z+c-b_m)}{s^2-\bar b_m(z+c)}\cdot\frac{s^2-
    \bar    b_m z}{s(z-b_m)}\right|,
    \end{split}
    \end{equation}
where $|z|=r$, $s>r+|c|$, and $\{a_j\}$ and $\{b_m\}$ are the
sequences of zeros and poles of $f$, respectively. By denoting
$\{q_k\}:=\{a_j\}\cup\{b_m\}$ and integrating
\eqref{integratethis} over the set $\{\xi\in[0,2\pi):
\left|\frac{f(re^{i\xi}+c)}{f(re^{i\xi})}\right|\geq 1\}$, it
follows that
    \begin{equation}\label{s1s2}
    m\left(r,\frac{f(z+c)}{f(z)}\right) \leq
     S_1(r)+S_2(r),
     \end{equation}
where
    \begin{equation}\label{S1prev}
    \begin{split}
    S_1(r) &= \int_0^{2\pi} \int_0^{2\pi}
    \left|\log|f(se^{i\theta})|\textrm{Re}\left(\frac{2c se^{i\theta}}{(se^{i\theta}-re^{i\varphi}-c)(se^{i\theta}-
    re^{i\varphi})}\right)\right|\,\frac{d\theta}{2\pi}\frac{d\varphi}{2\pi}
    \end{split}
    \end{equation}
and
    \begin{equation*}
    \begin{split}
    S_2(r) &= \sum_{|q_k|<s}
    \int_0^{2\pi}\log^{+}\left|1+\frac{c}{re^{i\varphi}-q_k}\right| \,\frac{d\varphi}{2\pi}
     +  \sum_{|q_k|<s}\int_0^{2\pi}\log^{+}\left|1-\frac{c}{re^{i\varphi}+c-q_k}\right|
    \,\frac{d\varphi}{2\pi}
    \\
    &\quad+ \sum_{|q_k|<s}  \int_0^{2\pi}\log^{+}\left|1+\frac{c}{re^{i\varphi}-\frac{s^2}{\bar q_k}
    }\right| \,\frac{d\varphi}{2\pi}
    +   \sum_{|q_k|<s} \int_0^{2\pi}\log^{+}\left|1-\frac{c}{re^{i\varphi}+c-\frac{s^2}{\bar q_k}
    }\right| \,\frac{d\varphi}{2\pi}.
    \end{split}
    \end{equation*}
By Fubini's theorem the order of integration in \eqref{S1prev} may
be changed, which results in
    \begin{equation}\label{S1early}
    \begin{split}
    S_1(r) &= \int_0^{2\pi} \left|\log|f(se^{i\theta})|\right| \int_0^{2\pi}
    \left|\textrm{Re}\left(\frac{2c se^{i\theta}}{(se^{i\theta}-re^{i\varphi}-c)(se^{i\theta}-
    re^{i\varphi})}\right)\right|\,\frac{d\varphi}{2\pi}\frac{d\theta}{2\pi}\\
    &\leq \frac{2|c|s}{(s-r-|c|)(s-r)^{1-\delta}} \int_0^{2\pi} \left|\log|f(se^{i\theta})|\right|  \int_0^{2\pi}
    \frac{1}{|se^{i\theta}-
    re^{i\varphi}|^\delta}\,\frac{d\varphi}{2\pi}\frac{d\theta}{2\pi}.
    \end{split}
    \end{equation}
By the change of variables $\varphi'=\theta-\varphi$, we have
    \begin{equation*}
    \begin{split}
    \int_0^{2\pi} \frac{1}{|se^{i\theta}-
    re^{i\varphi}|^\delta}\,\frac{d\varphi}{2\pi} &= \int_0^{2\pi}
    \frac{1}{|se^{i(\theta-\varphi)}-
    r|^\delta}\,\frac{d\varphi}{2\pi}= -\int_\theta^{\theta-2\pi}
    \frac{1}{|se^{i\varphi'}-r|^\delta}\,\frac{d\varphi'}{2\pi} \\ &= \int_0^{2\pi}
    \frac{1}{|se^{i\varphi'}- r|^\delta}\,\frac{d\varphi'}{2\pi}
    \leq \frac{2\pi}{s^\delta(1-\delta)}
    \end{split}
    \end{equation*}
(see, e.g., \cite[p. 89]{goldbergo:08} for the last inequality).
Hence \eqref{S1early} becomes
    \begin{equation}\label{S1}
    \begin{split}
    S_1(r) &\leq  \frac{4\pi|c|}{(1-\delta)(s-r-|c|)}\cdot
    \left(\frac{s}{s-r}\right)^{1-\delta}\left(m(s,f)+m\left(s,\frac{1}{f}\right)\right).
    \end{split}
    \end{equation}
Moreover, since
    \begin{equation*}
    \begin{split}
    \int_0^{2\pi}\log^+\left|1+\frac{c}{re^{i\varphi}-d}\right|\frac{\varphi}{2\pi}
    &\leq
    \frac{1}{\delta}\int_0^{2\pi}\log^+\left|1+\frac{c}{re^{i\varphi}-d}\right|^\delta\frac{\varphi}{2\pi}\\
    &\leq
    \frac{1}{\delta}\int_0^{2\pi}\left|\frac{c}{re^{i\varphi}-d}\right|^\delta\frac{\varphi}{2\pi}
    \leq\frac{2\pi|c|^\delta}{\delta(1-\delta)r^\delta}
    \end{split}
    \end{equation*}
for any $d\in\C$, it follows that
    \begin{equation}\label{S23}
    \begin{split}
    S_2(r) &\leq \frac{8\pi|c|^\delta}{\delta(1-\delta)r^\delta}
    \left(n(s,f)+n\left(s,\frac{1}{f}\right)\right).
    \end{split}
    \end{equation}
The assertion follows by combining the inequalities \eqref{s1s2},
\eqref{S1} and \eqref{S23}. \hfill $\Box$\par\vspace{2.5mm}

We will now extend Lemma~\ref{details1} to several complex variables.
The basic idea is to combine a method, which
Biancofiore and Stoll refer to as ``fiber integration''
\cite{biancofiores:81} (see also \cite{ye:96}) with Lemma
\ref{details1}. For the sake of brevity we adopt the notation
    \begin{equation*}
    \begin{split}
    m_f(r,\infty,0) &:= m(r,f)
    + m\left(r,\frac{1}{f}\right),\\
    n_f(r,\infty,0) &:= n_f(r,\infty)
    + n_f(r,0).\\
    \end{split}
    \end{equation*}

\begin{lemma}\label{prooflemma}
Let $f$ be a non-constant meromorphic function in $\C^n$ such that
$f(0)\not=0,\infty$, let $c=(c_1,\ldots,c_n)\in\C^n$, let
$0<\delta<1$, and denote $\tilde
c_j:=(0,\ldots,0,c_j,0,\ldots,0)$. Then
    \begin{equation*}
    \begin{split}
    &\int\limits_{\partial B_n(r)} \hspace{-2mm} \log^{+}\left|\frac{f(z+\tilde
    c_j)}{f(z)}\right|\sigma_n(z) \leq      \frac{8\pi|c_j|^\delta
    C}{\delta(1-\delta)}\left(\frac{R}{r}\right)^{2n-2}\frac{n_f(R,\infty,0)}{r^\delta} \\&\quad+
    \frac{4\pi|c_j|}{1-\delta}\left(\frac{R}{r}\right)^{2n-2}
     \left(\frac{R}{R-(r+|c_j|)}\right)\left(\frac{R}{R-r}\right)^{1-\delta}  \frac{m_f(R,\infty,0)}{\sqrt{R^2-r^2}}
     \end{split}
    \end{equation*}
for all $R>r+|c_j|>|c_j|$.
\end{lemma}

\begin{proof}
Let $r>0$, and let $h$ be a function on $\partial B_n(r)$ such
that $h\sigma_n$ is integrable over $\partial B_n(r)$. Then,
according to \cite[Lemma 3.1]{biancofiores:81},
    \begin{equation}\label{fiberint}
    \int\limits_{\partial B_n(r)}\hspace{-2mm} h(z)\sigma_n(z) =
    \frac{1}{r^{2n-2}}\int\limits_{\overline{B}_{n-1}(r)} \int\limits_{\partial
    B_1(p_r(w))}\hspace{-4mm}h(w,\zeta)\sigma_1(\zeta)\rho_{n-1}(w),
    \end{equation}
where $p_r(w)=\sqrt{r^2-|w|^2}$. Write $f_{[w]}(z)=f(w,z)$ for
$w\in\C^{n-1}$. By applying \eqref{fiberint} with
$h(z)=\log^{+}|f(z+\tilde c_j)/f(z)|$, we obtain
    \begin{equation}\label{estim1}
    \int\limits_{\partial B_n(r)} \hspace{-2mm} \log^{+}\left|\frac{f(z+\tilde c_j)}{f(z)}\right|\sigma_n(z) =
    \frac{1}{r^{2n-2}}\int\limits_{\overline{B}_{n-1}(r)} \int\limits_{\partial
    B_1(p_r(w))}\hspace{-4mm}\log^{+}\left|\frac{f_{[w]}(\zeta+
    c_j)}{f_{[w]}(\zeta)}\right|\sigma_1(\zeta)\rho_{n-1}(w).
    \end{equation}
Since $p_R(w)>p_r(w)+|c_j|$ whenever $R>r+|c_j|$, Lemma
\ref{details1}, applied with \eqref{estim1},  implies that
    \begin{equation}\label{I3}
    \begin{split}
    &\int\limits_{\partial B_n(r)} \hspace{-2mm} \log^{+}\left|\frac{f(z+\tilde
    c_j)}{f(z)}\right|\sigma_n(z)\\
    & \quad \leq
    \frac{1}{r^{2n-2}}\int\limits_{\overline{B}_{n-1}(r)}
    \Bigg(\frac{4\pi|c_j|}{(1-\delta)(p_R(w)-p_r(w)-|c|)}\cdot
    \left(\frac{p_R(w)}{p_R(w)-p_r(w)}\right)^{1-\delta} \\
    &\qquad\qquad\qquad\qquad\qquad\times  m_{f_{[w]}}(p_R(w),\infty,0)\Bigg)\rho_{n-1}(w)\\
     &\qquad+ \frac{1}{r^{2n-2}}\int\limits_{\overline{B}_{n-1}(r)}
    \frac{8\pi|c_j|^\delta}{\delta(1-\delta)p_r(w)^\delta}
    n_{f_{[w]}}(p_R(w),\infty,0)
    \rho_{n-1}(w)\\
    &\quad =:I_m+I_n
    \end{split}
    \end{equation}
for all $R>r+|c_j|$.

We will now proceed to estimate terms $I_m$ and $I_n$ separately,
starting with $I_m$. Since $p_r(w)+|c_j|\leq p_{r+|c_j|}(w)$ for
all $r>0$, and since $p_R(w)\geq \sqrt{R^2-r^2}$ and $p_r(w)\leq
r\cdot p_R(w)/R$ for all $w\in\overline{B}_{n-1}(r)$, it follows
that
    \begin{equation*}\label{pwestim}
     \frac{1}{p_R(w)-p_r(w)-|c_j|}\leq \frac{1}{\displaystyle p_R(w)\left(1-\frac{p_{r+|c_j|}(w)}{p_R(w)}\right)} \leq \frac{R}{(R-(r+|c_j|))\sqrt{R^2-r^2}}
    \end{equation*}
and
    \begin{equation*}
    \left(\frac{p_R(w)}{p_R(w)-p_r(w)}\right)^{1-\delta} \leq
    \left(\frac{R}{R-r}\right)^{1-\delta}.
    \end{equation*}
Therefore
    \begin{equation*}
    I_m \leq
     \left(\frac{R}{R-(r+|c_j|)}\right)\frac{4\pi|c_j|
     r^{2-2n}}{(1-\delta)\sqrt{R^2-r^2}} \left(\frac{R}{R-r}\right)^{1-\delta}
     \int\limits_{\overline{B}_{n-1}(R)} \hspace{-2mm}
      m_{f_{[w]}}(p_R(w),\infty,0)
    \rho_{n-1}(w).
    \end{equation*}
Since
    \begin{equation*}
    \frac{1}{R^{2n-2}} \int\limits_{\overline{B}_{n-1}(R)}\hspace{-2mm} m_{f_{[w]}}(p_R(w),\infty,0)
    \rho_{n-1}(w) = m_f(R,\infty,0)
    \end{equation*}
by equation \eqref{fiberint}, we finally have
    \begin{equation}\label{Im}
    I_m \leq \frac{4\pi|c_j|}{1-\delta}\left(\frac{R}{r}\right)^{2n-2}
     \left(\frac{R}{R-(r+|c_j|)}\right)\left(\frac{R}{R-r}\right)^{1-\delta}  \frac{m_f(R,\infty,0)}{\sqrt{R^2-r^2}}
    \end{equation}
for all $R>r+|c_j|$.

Consider now the term $I_{n}$. We may assume, without loss of
generality, that $\delta>1/4$. Then, denoting the integer part of
a real number $x$ by $[x]$, it follows that
$q(\delta):=[1/(1-\sqrt{\delta})]\geq 2$, and so H\"older's
inequality yields
    \begin{equation}\label{holder}
    \begin{split}
    \int\limits_{\overline{B}_{n-1}(r)}\hspace{-2mm}
    \frac{n_{f_{[w]}}(p_R(w),\infty,0)}{p_r(w)^\delta}\,
    \rho_{n-1}(w)&\leq
    \left(\int\limits_{\overline{B}_{n-1}(r)}\hspace{-2mm} n^{q(\delta)}_{f_{[w]}}(p_R(w),\infty,0)
    \rho_{n-1}(w)\right)^{\frac{1}{q(\delta)}}\\
    &\quad\times\left(\int\limits_{\overline{B}_{n-1}(r)}\hspace{-2mm} p_r(w)^{-\frac{\delta
    q(\delta)}{q(\delta)-1}}
    \rho_{n-1}(w)\right)^{\frac{q(\delta)-1}{q(\delta)}}.
    \end{split}
    \end{equation}
Since $0<\frac{\delta q(\delta)}{q(\delta)-1}<1$, it follows that
    \begin{equation}\label{tulo1}
    \int\limits_{\overline{B}_{n-1}(r)}\hspace{-2mm} p_r(w)^{-\frac{\delta
    q(\delta)}{q(\delta)-1}}
    \rho_{n-1}(w) \leq C r^{2n-2-\frac{\delta
    q(\delta)}{q(\delta)-1}}
    \end{equation}
where
    \begin{equation*}
    C=\int\limits_{\overline{B}_{n-1}(1)}\hspace{-2mm} \frac{1}{(1-\xi^2)^\frac{\delta
    q(\delta)}{2(q(\delta)-1)}} \rho_{n-1}(\xi).
    \end{equation*}
On the other hand, by \cite[Hilfssatz 7]{stoll:52} applied with a
weighted counting function $\widetilde n$ such that $\widetilde
n(r)=n^{q(\delta)}_f(R,\infty,0)$, it follows that
    \begin{equation}\label{tulo2}
    \begin{split}
    n^{q(\delta)}_f(R,\infty,0)&=\widetilde n(R)\\
    &\geq \frac{1}{R^{2n-2}}\int\limits_{\overline{B}_{n-1}(R)}\hspace{-2mm} \widetilde
    n_{f_{[w]}}(p_R(w))\rho_{n-1}(w)\\&\geq \frac{1}{R^{2n-2}}
     \int\limits_{\overline{B}_{n-1}(r)}\hspace{-2mm} n^{q(\delta)}_{f_{[w]}}(p_R(w),\infty,0)
    \rho_{n-1}(w).
    \end{split}
    \end{equation}
Finally, by \eqref{holder}, \eqref{tulo1} and \eqref{tulo2}, we
have
    \begin{equation}\label{In2}
    I_{n} \leq \frac{8\pi|c_j|^\delta
    C}{\delta(1-\delta)}\left(\frac{R}{r}\right)^{2n-2}
     \frac{n_f(R,\infty,0)}{r^\delta}.
    \end{equation}
The assertion of the lemma follows by combining the estimates
\eqref{I3}, \eqref{Im} and \eqref{In2}.
\end{proof}

\section{Proof of Theorem \ref{logdiff}}\label{proofsec}

Since
\begin{equation*}
n_f(r,\infty,0) \leq
\frac{R}{R-r}\left(N(R,f)+N\left(R,\frac{1}{f}\right)\right)
\end{equation*}
for all $R>r$, it follows by the first main theorem \eqref{1mt}
and Lemma \ref{prooflemma} that there exists a positive constant $K_1$, depending only
on $c_j$ and $\delta$, such that
    \begin{equation}\label{Tineq}
    \begin{split}
    &\int\limits_{\partial B_n(r)} \hspace{-2mm} \log^{+}\left|\frac{f(z+\tilde
    c_j)}{f(z)}\right|\sigma_n(z) \leq K_1K_2(r,R)\left(T(R,f)+\log\frac{1}{|f(0)|}\right)
     \end{split}
    \end{equation}
for all $R>r+|c_j|>|c_j|$, where
\begin{equation}
K_2(r,R)=\left(\frac{R}{r}\right)^{2n-2}
     \left(\frac{R}{R-(r+|c_j|)}\right)\left(\frac{1}{\sqrt{R^2-r^2}}\left(\frac{R}{R-r}\right)^{1-\delta} +
     \frac{1}{r^\delta}\right).
\end{equation}
Let $\xi(x)$ and $\phi(r)$ be positive, nondecreasing, continuous
functions defined for $e\leq x<\infty$ and $r_0\leq r<\infty$,
respectively, where $r_0$ is such that $T(r+|c|,f)\geq e$ for all
$r\geq r_0$.  Then by Hinkkanen's  Borel type growth lemma
\cite[Lemma 3]{hinkkanen:92} (see also \cite[Lemma
3.3.1]{cherryy:01})
    \begin{equation*}
    T\left(r+|c|+\frac{\phi(r+|c|)}{\xi(T(r+|c|,f))},f\right) \leq 2 T(r+|c|,f)
    \end{equation*}
for all $r$ outside of a set $E$ satisfying
    \begin{equation*}
    \int_{E\cap [r_0,s]} \frac{dr}{\phi(r)} \leq \frac{1}{\xi(e)}+ \frac{1}{\log
     2}\int_e^{T(s+|c|,f)}\frac{dx}{x\xi(x)}
    \end{equation*}
where $s<\infty$. Therefore, by choosing $\phi(r)=r$ and
$\xi(x)=(\log x)^{1+\tilde\varepsilon}$ with
$\tilde\varepsilon>0$, and defining
    \begin{equation}\label{alpha}
    R = r+|c_j|+\frac{r+|c_j|}{(\log T(r+|c_j|,f))^{1+\tilde\varepsilon}},
    \end{equation}
 we have
    \begin{equation}\label{estimate}
    T(R,f)=T\left(r+|c_j|+\frac{\phi(r+|c_j|)}{\xi(T(r+|c_j|,f))},f\right) \leq  2T(r+|c_j|,f)
    \end{equation}
for all $r$ outside of a set $E$ of finite logarithmic measure. By
substituting \eqref{alpha} and \eqref{estimate} into
\eqref{Tineq}, we obtain
\begin{equation}\label{estimeps}
\int\limits_{\partial B_n(r)} \hspace{-2mm} \log^{+}\left|\frac{f(z+\tilde
    c_j)}{f(z)}\right|\sigma_n(z)= o\left(\frac{T(r+|c_j|,f)(\log T(r+|c_j|,f))^{(1+\tilde\varepsilon)(\frac{5}{2}-\delta)}}{r^\delta}\right)
\end{equation}
where $r$ runs to infinity outside of an exceptional set of finite
logarithmic measure. (From now on $E\subset [1,+\infty)$ denotes a
set, which is not necessarily the same at each occurrence, but
which always has finite logarithmic measure.)

Since the hyper-order of $f$ is $\varsigma(f)=\varsigma$, we have
$\log T(r+|c_j|,f)\leq r^{\varsigma+\tilde\varepsilon}$ for all
$r$ sufficiently large. On the other hand, by \cite[Lemma
2.1]{halburdk:07PLMS} (see also \cite[Lemma~8.3]{halburdkt:09}), we
have $T(r+|c_j|,f)=T(r,f)+o(T(r,f))$ for all $r$ outside of an
exceptional set $E$ of finite logarithmic measure. Therefore,
\eqref{estimeps} yields
    \begin{equation}\label{estimnoeps}
   \int\limits_{\partial B_n(r)} \hspace{-2mm} \log^{+}\left|\frac{f(z+\tilde
    c_j)}{f(z)}\right|\sigma_n(z)= o\left(\frac{T(r,f)}{r^{\delta(1+\varsigma)-\frac{5+\varepsilon}{2}\varsigma}}\right),
    \end{equation}
where $\varepsilon>0$ is arbitrary small (and depends only on
$\tilde\varepsilon$), and $r\not\in E$.

In the general case any $c\in\C^n$ can be written as
$c=\sum_{j=0}^n \tilde c_j$ where $\tilde c_0:=0$. Therefore, by
\eqref{estimnoeps}, we have
    \begin{equation}\label{gen}
    \begin{split}
    \int\limits_{\partial B_n(r)} \hspace{-2mm}
    \log^{+}\left|\frac{f(z+c)}{f(z)}\right|\sigma_n(z) &=
    \int\limits_{\partial B_n(r)} \hspace{-2mm}
    \log^{+}\prod_{k=1}^n\left|\frac{f(z+\sum_{j=0}^k \tilde
    c_j)}{f(z + \sum_{j=0}^{k-1} \tilde
    c_j)}\right|\sigma_n(z)\\
    &\leq \sum_{k=1}^n
    \int\limits_{\partial B_n(r)} \hspace{-2mm}
    \log^{+}\left|\frac{f(z+\sum_{j=0}^k \tilde
    c_j)}{f(z + \sum_{j=0}^{k-1} \tilde
    c_j)}\right|\sigma_n(z)\\
    &= \sum_{k=1}^n o\left(\frac{T(r,f(z + \sum_{j=0}^{k-1} \tilde
    c_j))}{r^{\delta(1+\varsigma)-\frac{5+\varepsilon}{2}\varsigma}}\right)
    \end{split}
    \end{equation}
for all $r\not\in E$. On the other hand, by \cite[Lemma
2.1]{halburdk:07PLMS} (see also \cite[Lemma~8.3]{halburdkt:09}) it
follows that for any $s>0$ which does not depend on $r$ we have
$N(r+s,f)=N(r,f)+o(N(r,f))$, where $r\not\in E$. Hence, by
\eqref{estimnoeps}, we have
    \begin{equation}\label{Trcj}
    \begin{split}
    T(r,f(z+\tilde c_j)) &= m(r,f(z+\tilde c_j)) + N(r,f(z+\tilde
    c_j)) \\
    &\leq m\left(r,\frac{f(z+\tilde c_j)}{f(z)}\right)+ m(r,f) +
    N(r+|c_j|,f)\\
    &=T(r,f)+o(T(r,f))
    \end{split}
    \end{equation}
for all $r\not\in E$. Since  $c=\sum_{j=0}^n \tilde c_j$, it
follows by repeated application of \eqref{Trcj} that
    \begin{equation}\label{Trc}
    T(r,f(z+c)) = T(r,f) +o(T(r,f))
    \end{equation}
where $r\not\in E$. Relation \eqref{mainrel} follows by combining
\eqref{gen} and \eqref{Trc}, and by substituting $\delta=1-\varepsilon/(2+2\varsigma)$. \hfill $\Box$\par\vspace{2.5mm}

\section{Proof of Theorem \ref{2mthm}}\label{proof2nd}

The first main theorem yields
    \begin{equation}\label{2ndmaintheorem:3}
    \begin{split}
    \sum_{k=1}^pm\left(r,\frac{1}{f-a_k}\right)&=\sum_{k=1}^pT\left(r,\frac{1}{f-a_k}\right)
    -\sum_{k=1}^pN\left(r,\frac{1}{f-a_k}\right)\\
    &=pT(r,f)-N\left(r,\frac1{P(f)}\right)+O(1),
    \end{split}
    \end{equation}
where
\begin{equation*}
P(f)=\prod_{k=1}^p(f-a_k).
\end{equation*}
By partial fraction decomposition there exist constants
$\alpha_k\in\C$ such that
\begin{equation*}
\frac{1}{P(f)}=\sum_{k=1}^p\frac{\alpha_k}{f-a_k},
\end{equation*}
and so, since we have assumed that $f$ is finite at the origin and
$f(0)\not=a_j$ for $j=1,\ldots,q$, Corollary \ref{shiftcor} yields
\begin{equation*}
m\left(r,\frac{\Delta_c
f}{P(f)}\right)\leq\sum_{k=1}^pm\left(r,\frac{\Delta_c
f}{f-a_k}\right)+O(1)=o\left(\frac{T(r,f)}{r^\delta}\right)
\end{equation*}
for all $r$ outside of an exceptional set of finite logarithmic
measure. Therefore,
\begin{equation}
m\left(r,\frac{1}{P(f)}\right)=m\left(r,\frac{\Delta_cf}{P(f)}\frac{1}{\Delta_cf}\right)\leq
m\left(r,\frac{1}{\Delta_cf}\right) +
o\left(\frac{T(r,f)}{r^\delta}\right) \label{2ndmaintheorem:2}
\end{equation}
outside of an exceptional set. By applying Theorem~\ref{vm}, it
follows that $pT(r,f)=T(r,P(f))+O(1)$, and so by using the first
main theorem and \eqref{2ndmaintheorem:2},
Eq.~\eqref{2ndmaintheorem:3} becomes
\begin{eqnarray*}
\sum_{k=1}^pm\left(r,\frac{1}{f-a_k}\right)
&=&m\left(r,\frac{1}{P(f)}\right)+o\left(\frac{T(r,f)}{r^\delta}\right)\\
&\leq& m\left(r,\frac{1}{\Delta_cf}\right)+o\left(\frac{T(r,f)}{r^\delta}\right)\\
&=&T(r,\Delta_cf)-N\left(r,\frac1{\Delta_cf}\right)+o\left(\frac{T(r,f)}{r^\delta}\right),
\end{eqnarray*}
where $r$ runs to infinity outside of an exceptional set of finite
logarithmic measure. Therefore,
\begin{eqnarray*}
m(r,f) + \sum_{k=1}^pm\left(r,\frac{1}{f-a_k}\right)&\leq&T(r,f)+N(r,\Delta_cf)+m(r,\Delta_cf)\\
&&-N\left(r,\frac{1}{\Delta_cf}\right)-N(r,f)+o\left(\frac{T(r,f)}{r^\delta}\right)\\
\end{eqnarray*}
outside the exceptional set. Since
\begin{equation*}
m(r,\Delta_cf)=m\left(r,f\frac{\Delta_cf}{f}\right)\leq
m(r,f)+m\left(r,\frac{\Delta_cf}{f}\right)=m(r,f)+o\left(\frac{T(r,f)}{r^\delta}\right)
\end{equation*}
by Corollary \ref{shiftcor}, it follows that
\begin{eqnarray*}
m(r,f) + \sum_{k=1}^pm\left(r,\frac{1}{f-a_k}\right)&\leq&2T(r,f)+N(r,\Delta_cf)-N\left(r,\frac1{\Delta_cf}\right)\\
&&-2N(r,f)+o\left(\frac{T(r,f)}{r^\delta}\right)
\end{eqnarray*}
for all $r$ outside of an exceptional set of finite logarithmic
measure. \hfill $\Box$\par\vspace{2.5mm}

\section{Proof of Theorem \ref{picard}}\label{picardproof}

By composing $f$ with an appropriate M\"obius transformation, if
necessary, it may be assumed that $a_j\in\C$ and $f(0)\not=a_j$
for $j=1,2,3$. Consider the composition of $f$ with the function
$\tau(z)=z+c$. Since, by Theorem~\ref{logdiff},
    \begin{equation*}
    m(r,f\circ\tau)=m(r,f)+o(T(r,f)),
    \end{equation*}
and by Theorem~\ref{Trfc},
    \begin{equation*}
    T(r,f\circ\tau)=T(r,f)+o(T(r,f))
    \end{equation*}
for all $r$ outside of an exceptional set of finite logarithmic
measure, it follows that
    \begin{equation*}
    \begin{split}
     N(r,\Delta_c f)&\leq
     N(r,f\circ\tau)+N(r,f)\\
     &=T(r,f\circ\tau)+T(r,f)-m(r,f\circ\tau)-m(r,f)\\
     &=2T(r,f)-2m(r,f)+o(T(r,f))\\
     &= 2N(r,f)+o(T(r,f))
    \end{split}
    \end{equation*}
outside of an exceptional set. Therefore, by Theorem \ref{2mthm}
it follows that either
    \begin{equation}\label{contra}
    T(r,f)
    \leq  \sum_{k=1}^3 N\left(r,\frac{1}{f-a_k}\right)
    - N\left(r,\frac{1}{f\circ\tau-f}\right)
     + o(T(r,f))
    \end{equation}
for all $r$ outside of a small exceptional set, or
$f\circ\tau\equiv f$. Since by the assumption
$\tau(f^{-1}(\{a_j\}))\subset f^{-1}(\{a_j\})$ for $j=1,2,3$, it
follows that
    \begin{equation*}
    \sum_{k=1}^3 N\left(r,\frac{1}{f-a_k}\right)
    \leq N\left(r,\frac{1}{f\circ\tau-f}\right)
    \end{equation*}
and thus \eqref{contra} leads to a contradiction. Therefore
$f\equiv f\circ\tau$. \hfill $\Box$\par\vspace{2.5mm}



\def\cprime{$'$}
\providecommand{\bysame}{\leavevmode\hbox to3em{\hrulefill}\thinspace}
\providecommand{\MR}{\relax\ifhmode\unskip\space\fi MR }
\providecommand{\MRhref}[2]{%
  \href{http://www.ams.org/mathscinet-getitem?mr=#1}{#2}
}
\providecommand{\href}[2]{#2}

\end{document}